\documentclass[11pt,reqno]{amsart}
\usepackage{amssymb}
\newtheorem{thm}{Theorem}[section]
\newtheorem{prop}[thm]{Proposition}
\newtheorem{lemma}[thm]{Lemma}
\newtheorem{cor}[thm]{Corollary}

\def\cA{\mathcal{A}}
\def\bR{\mathbb{R}}
\def\bC{\mathbb{C}}
\def\bN{\mathbb{N}}
\def\bF{\mathbb{F}}
\def\bT{\mathbb{T}}
\def\ffi{\varphi}
\def\eps{\varepsilon}
\def\Tr{\mathrm{Tr}}
\def\U{\mathrm{U}}
\def\SU{\mathrm{SU}}
\def\cP{\mathcal{P}}
\def\bQ{\mathbf{Q}}
\def\bA{\mathbf{A}}
\def\bB{\mathbf{B}}
\def\Ric{\mathrm{Ric}}
\def\Hess{\mathrm{Hess}}
\def\im{\sqrt{-1}}

\begin{document}

\title[Free TCI for Non-commutative Multi-variables]
{Free Transportation Cost Inequalities for non-commutative multi-variables}
\author[F. Hiai]{Fumio Hiai$\,^{1,2}$}
\address{Graduate School of Information Sciences,
Tohoku University, Aoba-ku, Sendai 980-8579, Japan}
\author[Y. Ueda]{Yoshimichi Ueda$\,^{1,3}$}
\address{Graduate School of Mathematics,
Kyushu University, Fukuoka 810-8560, Japan}

\thanks{$^1\,$Supported in part by Japan Society for the Promotion of Science,
Japan-Hungary Joint Project.}
\thanks{$^2\,$Supported in part by Grant-in-Aid for Scientific Research
(C)14540198 and by Strategic Information and Communications R\&D Promotion
Scheme of MPHPT}
\thanks{$^3\,$Supported in part by Grant-in-Aid for Young
Scientists (B)14740118.}

\maketitle

\begin{abstract}
We prove the free analogue of the transportation cost inequality for tracial
distributions of non-commutative self-adjoint (also unitary) multi-variables
based on random matrix approximation procedure.
\end{abstract}

\section*{Introduction}

The transportation cost inequality (TCI) gives an upper bound for the
quadratic Wasserstein distance by the square roof of the relative entropy.
For probability measures on a Polish space $\mathcal{X}$, the relative
entropy is $S(\mu,\nu)=\int_{\mathcal{X}}\log{d\mu\over d\nu}\,d\mu$ if
$\mu\ll\nu$ (otherwise, $S(\mu,\nu)=+\infty$) while the quadratic Wasserstein
distance is defined as
\begin{equation}\label{F-0.1}
W_2(\mu,\nu):=\inf_\pi\sqrt{\iint_{\!\!\!\mathcal{X}\times\mathcal{X}}
d(x,y)^2\,d\pi(x,y)},
\end{equation}
where $d(x,y)$ is the metric on $\mathcal{X}$ and $\pi$ runs over the
probability measures on $\mathcal{X}\times\mathcal{X}$ with marginals $\mu$
and $\nu$. In 1996, M.~Talagrand \cite{T} obtained the celebrated TCI
$W_2(\mu,\nu)\le\sqrt{2S(\mu,\nu)}$ for probability measures on $\bR^n$, where
$\nu$ is the standard Gaussian measure on $\bR^n$. Since then, the TCI has
been received a lot of attention. It was shown by F.~Otto and C.~Villani
\cite{OV} that, in the Riemannian manifold setting, the TCI follows from the
logarithmic Sobolev inequality (LSI) of D.~Bakry and M.~Emery \cite{BE}. The LSI
gives a lower bound for the relative Fisher information by the relative entropy,
which has played important roles in several contexts. Recent developments in
both LSI and TCI are found in \cite{L1, Vi} for example.

On the other hand, Ph.~Biane and D.~Voiculescu \cite{BV} proved the free
analogue of Talagrand's TCI for compactly supported measures on $\bR$, where
the relative entropy is replaced by its free analogue and the Gaussian measure
by the semicircular one. In \cite{HPU1,HPU2} we developed the random matrix
approximation method to obtain a slight generalization of Biane and
Voiculescu's free TCI as well as its counterpart on the circle $\bT$. The free
analogues of the LSI's on $\bR$ and on $\bT$ were also obtained in \cite{Bi}
and \cite{HPU1,HPU3} by the same method.

Recently, M.~Ledoux \cite{L2} used a similar random matrix technique to prove
the free analogue of the Brunn-Minkowski inequality for measures on $\bR$,
from which (together with the Hamilton-Jacobi approach) he gave short proofs
of the free TCI and LSI for measures on $\bR$. Furthermore, his approach was
shown in \cite{HP2} to be still applicable for getting the free TCI in
\cite{HPU2} for measures on $\bT$.

The free TCI's and LSI's so far are restricted to measures on $\bR$ or $\bT$
and are not truly non-commutative. However, Voiculescu's free entropy for
multi-variables was well developed in \cite{V2} (see also \cite{V-LMS} and
\cite[Chap.~6]{HP}) and the Wasserstein distance was also introduced in
\cite{BV} for multi-variables in the $C^*$-algebra setting; so we must
be in a good position to extend the free TCI to non-commutative
multi-variables. This is what we are going to do here. In fact, we will show
the truly non-commutative free TCI when the ``reference distribution" is
chosen to be that associated with freely independent (self-adjoint or unitary)
random variables. It of course includes the above-mentioned free analogue of
Talagrand's TCI. However, the present work is still in a very beginning in
this direction of the subject matter. For example, it is interesting to seek
for a non-commutative generalization of the above-mentioned Ledoux's approach,
which probably brings a new insight into free probability theory.     
 
In this paper, after preliminaries on the Wasserstein distance in \S1
following \cite{BV}, we obtain in \S2 the free TCI for non-commutative tracial
distributions of self-adjoint multi-variables with respect to a certain free
product distribution (see Theorem \ref{T-2.2}). In \S3 we present a sharper
TCI (see Theorem \ref{T-3.1}) by replacing the free entropy with another free
entropy-like quantity (introduced in \cite{H} from the viewpoint of statistical
mechanics) but tracial distributions are rather restricted. Furthermore, the
counterparts of these free TCI's in the unitary setting are sketched in \S4
without much details for proofs.

\section{Preliminaries}
\setcounter{equation}{0} 

\subsection{Notations} When $\cA$ is a unital $C^*$-algebra, $\cA^{sa}$
stands for the set of self-adjoint elements  of $\cA$, and we denote by
$S(\cA)$ the state space of $\cA$  and by $TS(\cA)$ the tracial state space
of $\cA$, i.e., the set  of all $\ffi\in S(\cA)$ such that
$\ffi(ab)=\ffi(ba)$, $a,b\in\cA$.  The universal free product $C^*$-algebra
of two copies of $\cA$ is denoted by $\cA\bigstar\cA$, and $\sigma_1$ and 
$\sigma_2$ stand for the canonical embedding maps of $\cA$ into the left and
right copies of $\cA$ in $\cA\bigstar\cA$, respectively. Moreover, the
universal free product $C^*$-algebra of $n$ copies  of $\cA$ is simply written
as $\cA^{\bigstar n}$. A pair $(\cA,\tau)$  with $\tau\in TS(\cA)$ is called
a tracial $C^*$-probability space, and when $\cA$ is a von Neumann algebra
and $\tau$ is a faithful normal tracial state it is called a tracial
$W^*$-probability space. 

The usual non-normalized trace on the $N\times N$ complex matrix algebra
$M_N(\bC)$ is denoted by $\mathrm{Tr}_N$, and $\|A\|_{HS}$ is the
Hilbert-Schmidt norm of $A\in M_N(\bC)$, i.e.,
$\|A\|_{HS}:=\sqrt{\Tr_N(A^*A)}$ while $\|A\|$ is the operator norm of $A$.
We denote by $M_N^{sa}$ the set of all self-adjoint $A\in M_N(\bC)$ and by
$\Lambda_N$ the Lebesgue measure on $M_N^{sa}$ with the obvious Euclidean
structure. As usual, $\U(N)$ and $\SU(N)$ are the unitary and special
unitary groups of order $N$. We denote by $\gamma_N^\U$ and
$\gamma_N^\SU$  the Haar probability measures on $\U(N)$ and $\SU(N)$,
respectively. We also denote by $\cP(\mathcal{X})$ the set of all Borel
probability measures on a Polish space $\mathcal{X}$.

\subsection{Non-commutative distributions} Slightly unlike the usual, we will
employ the scheme in \cite{H} to deal with ``non-commutative distributions."
Let us fix $n \in \mathbb{N}$ and $R > 0$. An underlying $C^*$-algebra
we adopt is the universal free product $C^*$-algebra $\mathcal{A}^{(n)}_R :=
C([-R,R])^{\bigstar n}$ with norm $\|\cdot\|_R$ and a canonical
set of self-adjoint generators $X_i(t) = t$ in the $i$th copy of $C([-R,R])$,
$1\le i\le n$. Each $\varphi \in S\bigl(\mathcal{A}^{(n)}_R\bigr)$ provides
a {\it distribution} or {\it law} of $X_1,\dots,X_n$ whose (non-commutative)
{\it moments} are given by $\varphi\left(X_{i_1}\cdots X_{i_m}\right)$'s.
Any distribution in the $C^*$-algebra setting can be indeed realized in this
way. More precisely, if $a_1,\dots,a_n$ are self-adjoint variables in a
$C^*$-probability space $(\cA,\ffi)$ with operator norm $\|a_i\|\le R$, then
one has a (unique) $\ast$-homomorphism $\Psi$ from $\cA_R^{(n)}$ into $\cA$
sending each $X_i$ to $a_i$ so that the distribution of $X_1,\dots,X_n$ under
$\ffi\circ\Psi\in S\bigl(\cA_R^{(n)}\bigr)$ coincides with that of
$a_1,\dots,a_n$ under $\ffi$. Our main objects in the paper are the
Wasserstein distance and the free entropy, which have been well developed
only in terms of tracial states. Thus, in what follows, we will restrict our
consideration only to tracial distributions, i.e., elements in
$TS\bigl(\cA_R^{(n)}\bigr)$.

The (microstates) {\it free entropy} $\chi$ introduced by Voiculescu \cite{V2}
is defined in our context for every
$\tau\in TS\bigl(\mathcal{A}^{(n)}_R\bigr)$ as follows: Let $\pi_{\tau}$ be
the GNS representation of $\mathcal{A}^{(n)}_R$ associated with $\tau$; then
we have the tracial $W^*$-probability space
$\Bigl(\pi_{\tau}\bigl(\cA_R^{(n)}\bigr)'',\widetilde{\tau}\Bigr)$ with the
normal extension $\widetilde{\tau}$ of $\tau$ together with the self-adjoint
variables $\pi_{\tau}(X_1),\dots,\pi_{\tau}(X_n)$. Then, the free entropy of
$\tau$ at our disposal is
\begin{equation*} 
\chi(\tau) := \chi(\pi_{\tau}(X_1),\dots,\pi_{\tau}(X_n))
= \chi_R(\pi_{\tau}(X_1),\dots,\pi_{\tau}(X_n)) 
\end{equation*} 
(see \cite{V2} and also \cite[6.3.6]{HP} for the latter equality). By
definition the free entropy $\chi(\tau)$ is determined only by the moments of
$\tau$ in $(X_1,\dots,X_n)$ independently of a particular choice of $R>0$.
(This is also the case for the Wasserstein distance as will be seen in
\S\S1.3.)

Here, let us introduce a certain class of non-commutative distributions
coming from so-called matrix integrals, which will play an important
role in the paper. For each $N \in \mathbb{N}$ and $A_1,\dots,A_n \in
M_N^{sa}$ with $\Vert A_i\Vert \leq R$ we have the ``non-commutative
functional calculus" 
\begin{equation*} 
h \in \mathcal{A}_R^{(n)} \mapsto h(A_1,\dots,A_n) \in M_N(\bC)  
\end{equation*}
that is the canonical $*$-homomorphism from $\mathcal{A}_R^{(n)}$ into
$M_N(\bC)$ sending each $X_i$ to $A_i$. Let $r_R$ be the retraction of
$\bR$ onto $[-R,R]$, i.e., 
\begin{equation*} 
r_R(t) := \begin{cases} -R & \text{if $t < -R$}, \\
\ t & \text{if $-R \leq t \leq R$}, \\
\ R & \text{if $t > R$}. \end{cases}
\end{equation*} 
The next lemma is quite easy to show from the obvious inequality
$|r_R(\alpha)-r_R(\beta)| \leq |\alpha-\beta|$ for $\alpha,\beta \in \bR$; so
we omit the proof.  

\begin{lemma}\label{L-1.1}
We have  $\Vert r_R(A) - r_R(B) \Vert_{HS} \leq \Vert A - B \Vert_{HS}$
for every $A, B \in M_N^{sa}$.
\end{lemma}

Hence, a usual approximation argument shows that the function
$(A_1,\dots,A_n) \mapsto\allowbreak h(r_R(A_1),\allowbreak\dots,r_R(A_n))$
is continuous on $(M_N^{sa})^n \cong \bR^{N^2n}$ with respect to the
Euclidean structure for each fixed $h \in \mathcal{A}_R^{(n)}$. Thus, each
probability measure $\lambda \in \mathcal{P}\bigl((M_N^{sa})^n\bigr)$ gives
rise to the tracial distribution $\widehat{\lambda}_R \in
TS\bigl(\mathcal{A}_R^{(n)}\bigr)$ defined by 
\begin{equation*} 
\widehat{\lambda}_R(h) := \int_{(M_N^{sa})^n}\frac{1}{N}
\mathrm{Tr}_N(h(r_R(A_1),\dots,r_R(A_n)))\,d\lambda(A_1,\dots,A_n),
\quad h \in \mathcal{A}_R^{(n)}.
\end{equation*}
We call this $\widehat{\lambda}_R$ the {\it random matrix distribution}
associated with $\lambda$. When the measure $\lambda$ is supported in
$(M_{N,R}^{sa})^n$ where $M_{N,R}^{sa}:=\{A\in M_N^{sa}:\|A\|\le R\}$, the
retraction $r_R$ is of course not needed in the above definition so that
$\widehat{\lambda}_R$ is simply defined by integrating over $(M_{N,R}^{sa})^n$. 

\subsection{Wasserstein distance} This part is from \cite{BV} with slight
modifications. Let $(a_1,\dots,a_n)$ and $(b_1,\dots,b_n)$ be $n$-tuples of
non-commutative variables in tracial $C^*$-probability spaces $(\cA_1,\tau_1)$
and $(\cA_2,\tau_2)$, respectively. Here, it may be emphasized that $a_i$'s
as well as $b_i$'s are not necessarily self-adjoint (even not normal). We
write $(a_1,\dots,a_n)\sim(b_1,\dots,b_n)$ if the $*$-distributions (or
$*$-moments) of $(a_1,\dots,a_n)$ and of $(b_1,\dots,b_n)$ are same, i.e.,
$$
\tau_1\bigl(a_{i_1}^{\eps_1}\cdots a_{i_m}^{\eps_m}\bigr)
=\tau_2\bigl(b_{i_1}^{\eps_1}\cdots b_{i_m}^{\eps_m}\bigr)
$$
for all $m\in\bN$, $i_1,\dots,i_m\in\{1,\dots,n\}$ and
$\eps_1,\dots,\eps_m\in\{1,*\}$. For $1\le p<\infty$, the {\it $p$-Wasserstein
distance} introduced in \cite{BV} is defined by
\begin{equation}\label{F-1.1}
W_p((a_1,\dots,a_n),(b_1,\dots,b_n))
:=\inf\left\{\Biggl(\sum_{i=1}^n\tau(|a_i'-b_i'|^p)\Biggr)^{1/p}\right\},
\end{equation}
where infimum is taken over all $2n$-tuples
$(a_1',\dots,a_n',b_1',\dots,b_n')$ in some tracial $C^*$-probability space
$(\cA,\tau)$ such that $(a_1',\dots,a_n')\sim(a_1,\dots,a_n)$ and
$(b_1',\dots,b_n')\sim(b_1,\dots,b_n)$. The definition itself says that the
quantity \eqref{F-1.1} depends only on the $*$-moments of $(a_1,\dots,a_n)$
and $(b_1,\dots,b_n)$.

Another definition of $W_p$ was also introduced in \cite{BV}, which is a bit
more tractable than the above. Let $\cA$ be a unital $C^*$-algebra with a
specified $n$-tuple $(a_1,\dots,a_n)$ of generators. For a given 
pair $\tau_1, \tau_2 \in TS(\cA)$ we define the set of (non-commutative 
tracial) {\it joining states} between $\tau_1$ and $\tau_2$ by   
\begin{equation*} 
TS(\cA\bigstar\cA \,|\, \tau_1,\tau_2) := 
\left\{ \tau \in TS(\cA\bigstar\cA) : \tau\circ\sigma_1 = \tau_1,
\,\tau\circ\sigma_2 = \tau_2 \right\}.
\end{equation*} 
For $1 \leq p < \infty$, the {\it $p$-Wasserstein distance} between
$\tau_1$ and $\tau_2$ is defined by 
\begin{equation}\label{F-1.2}
W_p(\tau_1,\tau_2)
:= \inf\left\{ \Biggl(\sum_{i=1}^n \tau\bigl(|\sigma_1(a_i)
- \sigma_2(a_i)|^p\bigr)\Biggr)^{1/p} :
\tau \in TS(\cA\bigstar\cA \,|\, \tau_1,\tau_2) \right\}.
\end{equation}
As remarked in \cite[\S\S1.2]{BV}, the two definitions \eqref{F-1.1} and
\eqref{F-1.2} give the same quantity in the following way.

\begin{prop}
For every $\tau_1,\tau_2\in TS(\cA)$, let $(a_1',\dots,a_n')$ and
$(a_1'',\dots,a_n'')$ be $(a_1,\dots,a_n)$ in $(\cA,\tau_1)$ and in
$(\cA,\tau_2)$, respectively. Then we have
$$
W_p(\tau_1,\tau_2)=W_p((a_1',\dots,a_n'),(a_1'',\dots,a_n'')).
$$
\end{prop}

The proof is easily done by manipulating appropriate GNS representations so
that we leave it to the reader. An important  consequence of the proposition
is that  $W_p(\tau_1,\tau_2)$ in \eqref{F-1.2} is independent of a particular
choice of $\cA$ with  a specified $n$-tuple $(a_1,\dots,a_n)$; namely, it is
determined  only by the $*$-moments of $\tau_1$ and $\tau_2$ in
$(a_1,\dots,a_n)$.

Basic properties of $W_p$ are in order. 
\begin{itemize}
\item[$1^\circ$] $W_p(\tau_1,\tau_2)$ is a metric on $TS(\cA)$
(see \cite[Theorem 1.3]{BV}).
\item[$2^\circ$] $W_p(\tau_1,\tau_2)$ is jointly lower semi-continuous in
$(\tau_1,\tau_2)\in TS(\cA)\times TS(\cA)$ in weak*-topology
(see \cite[Proposition 1.4]{BV}).
\item[$3^\circ$] $W_p(\tau_1,\tau_2)^p$ is jointly convex in
$(\tau_1,\tau_2)\in TS(\cA)\times TS(\cA)$. (This is easy to prove though not
included in \cite{BV}.)
\item[$4^\circ$] If $a_1,\dots,a_n$ are self-adjoint (or more generally
normal) and mutually commuting, then $W_p(\tau_1,\tau_2)$ coincides with the
usual $p$-Wasserstein distance $W_p(\mu_1,\mu_2)$ (see \eqref{F-0.1}), where
$\mu_1,\mu_2\in\cP(\bR^n)$ (or $\cP(\bC^n)$) are the spectral distribution
measures of the $n$-tuple $(a_1,\dots,a_n)$ constructed via the GNS
representations associated with $\tau_1,\tau_2$, respectively (see
\cite[Theorem 1.5]{BV}).
\end{itemize}

We will treat the (quadratic) $2$-Wasserstein distance $W_2(\tau_1,\tau_2)$
for tracial distributions of self-adjoint random variables in
\S2, \S3 and for those of unitary random variables in \S4. In the
self-adjoint case, we will take the universal $\mathcal{A}_R^{(n)}$ with the
specified self-adjoint generators $X_1,\dots,X_n$. This is indeed
universal in the sense that when $a_1,\dots,a_n$ are self-adjoint in any
$\cA$, any tracial distribution of $(a_1,\dots,a_n)$ can be realized via some
$\tau\in TS\bigl(\cA_R^{(n)}\bigr)$ as long as $R\ge\|a_i\|$, $1\le i\le n$
(see \S\S1.2).

In this subsection, we provide an inequality between the
free and usual $2$-Wasserstein distances for random matrix distributions
introduced in \S\S1.2, which will be  one of the keys in our later
discussions. The inequality corresponds to that in \cite[Lemmas 2.6 and
2.8]{HPU2}; however the argument here  is simpler than there because
we do not (indeed cannot) treat the ``eigenvalue distributions."
 
\begin{lemma}\label{L-1.3}
For every pair $\lambda_1,\lambda_2 \in \mathcal{P}\bigl((M_N^{sa})^n\bigr)$
and every $R>0$, let $\widehat{\lambda}_{1,R},\widehat{\lambda}_{2,R} \in
TS\bigl(\mathcal{A}_R^{(n)}\bigr)$ be the corresponding random matrix
distributions. Then we have
\begin{equation*} 
W_2\bigl(\widehat{\lambda}_{1,R},\widehat{\lambda}_{2,R}\bigr)
\leq \frac{1}{\sqrt{N}}W_2(\lambda_1,\lambda_2), 
\end{equation*} 
where $W_2(\lambda_1,\lambda_2)$ is the usual $2$-Wasserstein distance between
$\lambda_1, \lambda_2$ defined by
\begin{equation*} 
\inf_\pi\sqrt{\iint_{\!\!\!(M_N^{sa})^n\times(M_N^{sa})^n}
\sum_{i=1}^n \Vert A_i - B_i\Vert_{HS}^2\,d\pi} 
\end{equation*}
over the joining measures $\pi$ on $(M_N^{sa})^n\times(M_N^{sa})^n$ of
$\lambda_1, \lambda_2$, i.e., measures whose marginals are $\lambda_1$,
$\lambda_2$ {\rm(}see also \eqref{F-0.1}{\rm)}.
\end{lemma}     
\begin{proof} For each $n$-tuple $\bA=(A_1,\dots,A_n) \in (M_N^{sa})^n$ one has
the $*$-homomorphism
\begin{equation*} 
\Psi_R^\bA :  h \in \mathcal{A}_R^{(n)} \mapsto
h(r_R(A_1),\dots,r_R(A_n)) \in M_N(\mathbf{C})
\end{equation*} 
sending $X_i$ to $r_R(A_i)$ (see \S\S1.2), and moreover for each
$\bA,\bB\in(M_N^{sa})^n$ there is a unique
$*$-homomorphism 
\begin{equation*} 
\Psi_R^{\bA,\bB} := \psi_R^\bA\bigstar\psi_R^\bB
: \mathcal{A}_R^{(n)}\bigstar\mathcal{A}_R^{(n)} \rightarrow M_N(\mathbf{C})
\end{equation*} 
determined by 
\begin{align*} 
\Psi_R^{\bA,\bB}\circ\sigma_1 = \Psi_R^\bA, \quad
\Psi_R^{\bA,\bB}\circ\sigma_2 = \Psi_R^\bB.
\end{align*}
As in \S\S1.2, the function 
$(\bA,\bB) \mapsto \Psi_R^{\bA,\bB}(h)$
is continuous with respect to the Hilbert-Schmidt norms for each fixed
$h \in \mathcal{A}_R^{(n)}\bigstar\mathcal{A}_R^{(n)}$; hence every
joining measure $\pi$ of $\lambda_1, \lambda_2$ gives rise to the tracial
distribution $\widehat{\pi}_R \in
TS\bigl(\mathcal{A}_R^{(n)}\bigstar\mathcal{A}_R^{(n)}\bigr)$ defined by 
\begin{equation*} 
\widehat{\pi}_R(h) := \iint_{\!\!\!(M_N^{sa})^n\times(M_N^{sa})^n}
\frac{1}{N}\mathrm{Tr}_N\bigl(\Psi_R^{\bA,\bB}(h)\bigr)\,d\pi(\bA,\bB),
\quad h \in \mathcal{A}_R^{(n)}\bigstar\mathcal{A}_R^{(n)},
\end{equation*} 
which satisfies
\begin{equation*} 
\widehat{\pi}_R\circ\sigma_1 = \widehat{\lambda}_{1,R}, \quad
\widehat{\pi}_R\circ\sigma_2 = \widehat{\lambda}_{2,R}.
\end{equation*} 
Therefore, we have $\widehat{\pi}_R \in
TS\bigl(\mathcal{A}_R^{(n)}\bigstar\mathcal{A}_R^{(n)}
\big|\widehat{\lambda}_{1,R},\widehat{\lambda}_{2,R}\bigr)$ so that
\allowdisplaybreaks{
\begin{eqnarray*} 
W_2\bigl(\widehat{\lambda}_{1,R},\widehat{\lambda}_{2,R}\bigr)^2
&\leq& \widehat{\pi}_R\Biggl(\sum_{i=1}^n
(\sigma_1(X_i) - \sigma_2(X_i))^2\Biggr) \\
&=& \iint \sum_{i=1}^n \frac{1}{N}\Tr_N\bigl(\Psi_R^{\bA,\bB}
\bigl((\sigma_1(X_i) - \sigma_2(X_i))^2\bigr)\bigr)\,d\pi(\bA,\bB) \\
&=&\iint \sum_{i=1}^n\frac{1}{N}\Tr_N\bigl((r_R(A_i) - r_R(B_i))^2\bigr)
\,d\pi(\bA,\bB) \\
&=& \iint \sum_{i=1}^n\frac{1}{N} \Vert r_R(A_i) - r_R(B_i)\Vert_{HS}^2 
\,d\pi(\bA,\bB) \\
&\leq& \frac{1}{N} \iint \sum_{i=1}^n \Vert A_i - B_i\Vert_{HS}^2 
\,d\pi(\bA,\bB),
\end{eqnarray*}
}where the latter inequality is due to Lemma \ref{L-1.1}. Hence,
the desired inequality follows by taking the infimum of the last integral
over the joining measures $\pi$ of $\lambda_1,\lambda_2$.
\end{proof}

Finally, we remark that the $2$-Wasserstein distance is sometimes
defined with the cost function of the form
${1\over2}\times({\rm distance})^2$. In fact, in \cite{HP2,HPU1,HPU2} we adopted
the definition with a ${1\over2}$-multiple constant so that the bounds of
TCI's there and in the present paper are $2$ times different.

\section{Free TCI for $\chi$}
\setcounter{equation}{0}

We will obtain the free TCI for non-commutative tracial
distributions with respect to the distribution of freely independent random
variables, including a natural free analogue of celebrated Talagrand's TCI
\cite{T} with respect to the standard Gaussian measure on $\bR^n$. 
 
Let $\bQ = (Q_1,\dots,Q_n)$ be an $n$-tuple of real-valued continuous
functions on $\bR$ with 
\begin{equation}\label{F-2.1} 
\lim_{|x|\rightarrow\infty}\exp(-\varepsilon Q_i(x))=0 \quad
\text{for every $\varepsilon > 0$}. 
\end{equation} 
Then, for each $Q_i$ we define the $N\times N$ self-adjoint random matrix
$\lambda_N(Q_i) \in \mathcal{P}(M_N^{sa})$ by 
\begin{equation*} 
d\lambda_N(Q_i)(A) :=
\frac{1}{Z_N(Q_i)}\exp\bigl(-N\mathrm{Tr}_N(Q_i(A))\bigr)\,d\Lambda_N(A)
\end{equation*} 
with a normalization constant $Z_N(Q_i)$, whose mean eigenvalue distribution
on $\bR$ is denoted by $\widehat{\lambda}_N(Q_i)$. With $Q:=Q_i$, a
fundamental result in the theory of weighted potentials (see \cite[I.1.3]{ST})
tells us that the functional
$$
-\Sigma(\mu)+\mu(Q)
:=-\iint_{\!\!\!\bR^2}\log|x-y|\,d\mu(x)\,d\mu(y)+\int_\bR Q(x)\,d\mu(x),
\quad \mu\in\cP(\bR),
$$
has a unique minimizer $\mu_Q$ which is compactly supported and called the
{\it equilibrium measure} associated with $Q$. For example, when
$Q(x)=x^2/2$, $\mu_{Q}$ is the $(0,1)$-semicircular distribution
$d\gamma_{0,2}(x) := \frac{1}{2\pi}\sqrt{4-x^2}dx$ supported on
$[-2,2]$. Furthermore, the large deviation principle for self-adjoint random
matrices (see \cite{BG}, \cite[5.4.3]{HP}) shows that
$\widehat{\lambda}_N(Q)$ weakly converges to the equilibrium measure $\mu_Q$.
Let $R_0>0$ be the smallest such that all $\mu_{Q_i}$'s are supported in
$[-R_0,R_0]$; for example,
$R_0 = 2$ when $Q_i(x)=x^2/2$ for all $i$. We then notice that
\begin{equation}\label{F-2.2} 
\lambda_N(Q_i)(M_{N,R_0}^{sa}) =
\widehat{\lambda}_N(Q_i)([-R_0,R_0]) \longrightarrow
\mu_{Q_i}([-R_0,R_0]) = 1  
\end{equation}
as $N \rightarrow \infty$ for every $1\le i\le n$. Let us consider the product
measure    
$$ 
\lambda_N(\bQ) := \bigotimes_{i=1}^n \lambda_N(Q_i)
\in \mathcal{P}\bigl((M_N^{sa})^n\bigr),
$$
that is,
$$
d\lambda_N(\bQ)(A_1,\dots,A_n)
= \frac{1}{Z_N(\bQ)}
\exp\Biggl(-N\sum_{i=1}^n\mathrm{Tr}_N (Q_i(A_i))\Biggr)
d\Lambda_N^{\otimes n}(A_1,\dots,A_n)
$$
with $Z_N(\bQ):=\prod_{i=1}^n Z_N(Q_i)$, and
$\widehat{\lambda}_{N,R}(\bQ) \in TS\bigl(\mathcal{A}_R^{(n)}\bigr)$ denotes
the random matrix distribution associated with $\lambda_N(\bQ)$ (see \S\S1.2).
Furthermore, when $R \geq R_0$, define the tracial distribution
$\tau_\bQ \in TS\bigl(\mathcal{A}_R^{(n)}\bigr)$ to be the free product
state $\bigstar_{i=1}^n \mu_{Q_i}$ on
$\mathcal{A}_R^{(n)} = C([-R,R])^{\bigstar n}$, where each $\mu_{Q_i}$ is
meant a state on $C([-R,R])$ defined by integration. (Note that
the moments of $\tau_\bQ$ is independent of a choice of $R\geq R_0$.) 

We begin by restating the so-called asymptotic freeness due to Voiculescu
\cite{Vo} in our situation.

\begin{lemma}\label{L-2.1}
Whenever $R\geq R_0$ we have 
\begin{equation*} 
\lim_{N\rightarrow\infty}\widehat{\lambda}_{N,R}(\bQ) =
\tau_\bQ \quad \text{weakly*}. 
\end{equation*} 
\end{lemma}
 
\begin{proof} Fix an arbitrary $R \geq R_0$. By \eqref{F-2.2} we get 
\begin{equation*} 
\lambda_N(\bQ)\bigl((M_{N,R}^{sa})^n\bigr)
= \prod_{i=1}^n\widehat{\lambda}_N(Q_i)([-R,R])
\longrightarrow 1 \quad \text{as $N\rightarrow\infty$}. 
\end{equation*} 
For any non-commutative polynomial $p$ in $X_1,\dots,X_n$
($\in \mathcal{A}_R^{(n)}$) we have 
\allowdisplaybreaks{ 
\begin{align*}
\widehat{\lambda}_{N,R}(\bQ)(p)
&= \int \frac{1}{N}\mathrm{Tr}_N(p(r_R(A_1),\dots,r_R(A_n)))
\,d\lambda_N(\bQ)(A_1,\dots,A_n) \\
&= \int_{(M_{N,R}^{sa})^n}
\frac{1}{N}\mathrm{Tr}_N(p(A_1,\dots,A_n))
\,d\lambda_N(\bQ)(A_1,\dots,A_n) \\
&\phantom{aaa}+ \int_{(M_N^{sa})^n \setminus(M_{N,R}^{sa})^n}
\frac{1}{N}\mathrm{Tr}_N(p(r_R(A_1),\dots,r_R(A_n)))
\,d\lambda_N(\bQ)(A_1,\dots,A_n).   
\end{align*}
}Since 
\allowdisplaybreaks{
\begin{align*}  
&\left|\int_{(M_N^{sa})^n \setminus(M_{N,R}^{sa})^n}
\frac{1}{N}\mathrm{Tr}_N(p(r_R(A_1),\dots,r_R(A_n)))
\,d\lambda_N(\bQ)(A_1,\dots,A_n)\right| \\
&\phantom{aaaaaa} \leq \Vert p \Vert_R
\left(1- \lambda_N(\bQ)\bigl((M_{N,R}^{sa})^n\bigr)\right)
\longrightarrow 0 \quad \text{as $N \rightarrow \infty$},
\end{align*}
}the desired assertion follows from the naturally expected fact that 
\allowdisplaybreaks{
\begin{align*} 
\tau_\bQ(p) 
&= \lim_{N\rightarrow\infty} \int_{(M_{N,R}^{sa})^n}
\frac{1}{N}\mathrm{Tr}_N(p(A_1,\dots,A_n))
\,d\lambda_{N,R}(\bQ)(A_1,\dots,A_n) \\
&= \lim_{N\rightarrow\infty} \int_{(M_{N,R}^{sa})^n}
\frac{1}{N}\mathrm{Tr}_N(p(A_1,\dots,A_n))
\,d\lambda_N(\bQ)(A_1,\dots,A_n),  
\end{align*}
}where   
\allowdisplaybreaks{
\begin{gather*} 
\lambda_{N,R}(\bQ) := 
\frac{1}{\lambda_N(\bQ)\bigl((M_{N,R}^{sa})^n\bigr)} 
\lambda_N(\bQ)\Big|_{(M_{N,R}^{sa})^n}
= \bigotimes_{i=1}^n \lambda_{N,R}(Q_i), \\
\lambda_{N,R}(Q_i) := 
\frac{1}{\lambda_N(Q_i)\bigl(M_{N,R}^{sa}\bigr)}
\lambda_N(Q_i)\Big|_{M_{N,R}^{sa}}.
\end{gather*}
}Indeed, this is a simple consequence of an asymptotic freeness result in
\cite[4.3.5]{HP}, slightly generalizing Voiculescu's original in \cite{Vo}
to the setup in almost sure sense as well as to general unitarily invariant
self-adjoint random matrices. Also, one should note that $\lambda_{N,R}(Q_i)$
still weakly converges to $\mu_{Q_i}$ thanks to \eqref{F-2.2}. 
\end{proof}   

We are now in a position to state the main result of this section. 

\begin{thm}\label{T-2.2}
Assume that there exists a constant $\rho > 0$ such that all
$Q_i(x) - \frac{\rho}{2}x^2$, $1\le i\le n$, are convex on $\mathbf{R}$
{\rm(}so that the condition \eqref{F-2.1} automatically holds{\rm)}. Assume
$R \geq R_0$ with $R_0$ given above. Then we have   
\begin{equation*} 
W_2(\tau,\tau_\bQ) \leq
\sqrt{\frac{2}{\rho}\Biggl(-\chi(\tau)
+\tau\Biggl(\sum_{i=1}^n Q_i(X_i)\Biggr) + B_\bQ\Biggr)}
\end{equation*} 
for every $\tau \in TS\bigl(\mathcal{A}_R^{(n)}\bigr)$, where 
\begin{equation}\label{F-2.3}
B_\bQ := \lim_{N\rightarrow\infty}
\Biggl(\frac{1}{N^2}\sum_{i=1}^n\log Z_N(Q_i)
+\frac{n}{2}\log N\Biggr). 
\end{equation}
\end{thm}  

Since the equilibrium measure with respect to $Q(x)=x^2/2$ is the
$(0,1)$-semicircular distribution $\gamma_{0,2}$ and
\begin{equation*} 
\lim_{N\rightarrow\infty}\biggl(\frac{1}{N^2}\log Z_N(Q)
+\frac{1}{2}\log N\biggr) = \frac{1}{2}\log 2\pi
\end{equation*}
(see e.g.~\cite[4.4.6 and pp.~185--186]{HP}),
the above theorem includes the free analogue of Talagrand's TCI as follows.

\begin{cor}\label{C-2.3} 
If $R \geq 2$ and $\gamma_{0,2}^{\bigstar n} \in
TS\bigl(\mathcal{A}_R^{(n)}\bigr)$ is the {\rm (}non-commutative{\rm )}
distribution of a standard semicircular system, then 
\begin{equation*} 
W_2\Bigl(\tau,\gamma_{0,2}^{\bigstar n}\Bigr)
\leq \sqrt{2\Biggl(-\chi(\tau)
+ \tau\Biggl(\frac{1}{2}\sum_{i=1}^n X_i^2\Biggr)
+ \frac{n}{2}\log2\pi\Biggr) }
\end{equation*}
for every $\tau \in TS\bigl(\mathcal{A}_R^{(n)}\bigr)$. 
\end{cor} 

To prove Theorem \ref{T-2.2}, we need the following: 

\begin{lemma}\label{L-2.4}
Assume the same assumption for $Q_i$'s with a constant $\rho > 0$. 
Then, for every $N \in \bN$ and every $\lambda \in
\mathcal{P}\bigl((M_N^{sa})^n\bigr)$, we have 
\begin{equation*} 
W_2(\lambda, \lambda_N(\bQ)) \leq \sqrt{\frac{2}{\rho N}
S(\lambda, \lambda_N(\bQ))}, 
\end{equation*}
where $S(\lambda, \lambda_N(\bQ))$ is the relative entropy of $\lambda$ with
respect to $\lambda_N(\bQ)$. 
\end{lemma}
 
\begin{proof} Since all $Q_i(x) - \frac{\rho}{2}x^2$ are convex on
$\mathbf{R}$, so is 
\begin{equation*} 
\left(A_1,\dots,A_n\right) \in (M_N^{sa})^n
\mapsto  N \sum_{i=1}^n \mathrm{Tr}_N
\Bigl(Q_i(A_i)-\frac{\rho}{2}A_i^2 \Bigr).
\end{equation*}
(This is the reason why the multiple constant $1/N$ appears, see 
\cite[p.~212]{HPU2}.) Hence, the TCI for measures on Euclidean
spaces (see \cite[Theorem 6.5]{L1}) slightly generalizing Talagrand's
original implies the desired inequality with regarding
$(M_N^{sa})^n$ as $\bR^{N^2 n}$.
\end{proof}

\begin{proof}[Proof of Theorem \ref{T-2.2}]
First, note that the existence of the limit in \eqref{F-2.3} is in
\cite[5.4.3]{HP}. When $\chi(\tau) = -\infty$
nothing has to be done so that let us assume
$\chi(\tau) > -\infty$. Recall that 
\allowdisplaybreaks{
\begin{align*} 
\chi(\tau)
&= \chi_R(\pi_\tau(X_1),\dots,\pi_\tau(X_n)) \\
&= \lim_{\substack{m\rightarrow\infty \\
\varepsilon\searrow0}}\limsup_{N\rightarrow\infty}\left(
\frac{1}{N^2}\log\Lambda_N^{\otimes n}\bigl(
\Gamma_R(\pi_\tau(X_1),\dots,\pi_\tau(X_n);N,m,\varepsilon)\bigr)
+\frac{n}{2}\log N\right),
\end{align*}
}where $\Gamma_R(\pi_\tau(X_1),\dots,\pi_\tau(X_n);N,m,\varepsilon)$ is the
set of all $n$-tuples $(A_1,\dots,A_n) \in (M_{N,R}^{sa})^n$ such
that 
\begin{equation*} 
\left|\frac{1}{N}\mathrm{Tr}_N(A_{i_1}\cdots A_{i_r}) 
- \tau(X_{i_1}\cdots X_{i_r})\right| = 
\left|\frac{1}{N}\mathrm{Tr}_N(A_{i_1}\cdots A_{i_r}) 
- \widetilde{\tau}(\pi_\tau(X_{i_1})\cdots\pi_\tau(X_{i_r}))\right| 
< \varepsilon
\end{equation*} 
for all possible $i_1,\dots,i_r$ with $1 \leq r \leq m$. A suitable
subsequence $N(1)<N(2)<\cdots$ can be chosen in such a way that letting
\begin{equation*} 
\Gamma_R(\tau;k) :=
\Gamma_R\bigl(\pi_\tau(X_1),\dots,\pi_\tau(X_n);N(k),k,1/k\bigr) 
\end{equation*}
we get 
\begin{equation}\label{F-2.4}
\chi(\tau) = \lim_{k\rightarrow\infty} 
\biggl(\frac{1}{N(k)^2} \log\Lambda_{N(k)}^{\otimes n}(\Gamma_R(\tau;k)) +
\frac{n}{2}\log N(k)\biggr).
\end{equation}
Looking at this, we introduce the random matrix distribution
$\widehat{\lambda}_{N(k),R} \in TS\bigl(\mathcal{A}_R^{(n)}\bigr)$ associated
with the probability measure
\begin{equation*} 
\lambda_{N(k)} := \frac{1}{\Lambda_{N(k)}^{\otimes n}(\Gamma_R(\tau;k))}
\Lambda_{N(k)}^{\otimes n}\big|_{\Gamma_R(\tau;k)} \in
\mathcal{P}\bigl((M_{N(k),R}^{sa})^n\bigr).
\end{equation*}
Let $h$ be an arbitrary monomial $X_{i_1}\cdots X_{i_r} \in
\mathcal{A}_R^{(n)}$. As long as $r \leq k$, we get 
\begin{eqnarray*}
&&\left|\frac{1}{N(k)}\mathrm{Tr}_{N(k)}(h(A_1,\dots,A_n))-\tau(h)\right| \\
&&\qquad = \left|\frac{1}{N(k)}\mathrm{Tr}_{N(k)}(A_{i_1}\cdots A_{i_r})
-\tau(X_{i_1}\cdots X_{i_r})\right| < \frac{1}{k}
\end{eqnarray*} 
for all $(A_1,\dots,A_n) \in \Gamma_R(\tau;k)$, and hence
\begin{eqnarray*} 
&&\left|\widehat{\lambda}_{N(k),R}(h) - \tau(h)\right| \\
&&\qquad \leq 
\int_{\Gamma_R(\tau;k)}\left|\frac{1}{N(k)}\mathrm{Tr}_{N(k)}
(h(A_1,\dots,A_n))-\tau(h)\right| d\lambda_k(A_1,\dots,A_n) < \frac{1}{k}. 
\end{eqnarray*} 
This shows that $\widehat{\lambda}_{N(k),R}(h) \longrightarrow \tau(h)$ as
$k\rightarrow\infty$ for all monomials $h$ so that we get 
\begin{equation}\label{F-2.5}
\lim_{k\rightarrow\infty}\widehat{\lambda}_{N(k),R} = \tau
\quad \text{weakly*}.
\end{equation} 

By Lemmas \ref{L-1.3} and \ref{L-2.4} we have
\allowdisplaybreaks{
\begin{align*} 
W_2&\bigl(\widehat{\lambda}_{N(k),R},\widehat{\lambda}_{N(k),R}
(\bQ)\bigr)^2 \\
&\leq \frac{1}{N(k)} W_2(\lambda_{N(k)},\lambda_{N(k)}(\bQ))^2 \\
&\leq \frac{2}{\rho N(k)^2} S(\lambda_{N(k)},\lambda_{N(k)}(\bQ)) \\
&= \frac{2}{\rho N(k)^2} \int_{(M_{N(k),R}^{sa})^n}
\log \frac{d\lambda_{N(k)}}{d\lambda_{N(k)}(\bQ)}(A_1,\dots,A_n)
\,d\lambda_{N(k)}(A_1,\dots,A_n) \\
&= \frac{2}{\rho N(k)^2} 
\int_{(M_{N(k),R}^{sa})^n} 
\Biggl(-\log\Lambda_{N(k)}^{\otimes n}(\Gamma_R(\tau;k)) \\
&\phantom{aaaaaaa}
+  N(k)\sum_{i=1}^n\mathrm{Tr}_{N(k)}(Q_i(A_i))
+ \sum_{i=1}^n \log Z_{N(k)}(Q_i)\Biggr)
d\lambda_{N(k)}(A_1,\dots,A_n) \\ 
&= \frac{2}{\rho} \Biggl\{
- \biggl(\frac{1}{N(k)^2}\log\Lambda_{N(k)}^{\otimes n}(\Gamma_R(\tau;k))
+ \frac{n}{2}\log N(k)\biggr) \\
&\phantom{aaaaaa} +
\widehat{\lambda}_{N(k),R}\Biggl(\sum_{i=1}^n Q_i(A_i)\Biggr)
+ \Biggl(\frac{1}{N(k)^2}\sum_{i=1}^n\log Z_{N(k)}(Q_i)
+  \frac{n}{2}\log N(k)\Biggr)\Bigg\}.
\end{align*}
}The above last formula converges to
$$
\frac{2}{\rho} \Biggl(-\chi(\tau)
+ \tau\Biggl(\sum_{i=1}^n Q_i(A_i)\Biggr) + B_\bQ\Biggr)
$$
as $k\rightarrow\infty$ thanks to \eqref{F-2.4} and \eqref{F-2.5}.
On the other hand, by the joint lower semi-continuity of $W_2$ (see $2^\circ$
in \S\S1.3), Lemma 2.1 and \eqref{F-2.5}, we have
\begin{equation*} 
W_2(\tau,\tau_\bQ) \leq \liminf_{k\rightarrow\infty}W_2
\bigl(\widehat{\lambda}_{N(k),R},\widehat{\lambda}_{N(k),R}(\bQ)\bigr),
\end{equation*} 
completing the proof.   
\end{proof}

\section{Free TCI for $\eta$}
\setcounter{equation}{0}

We first recall the free pressure $\pi_R$ and the free entropy-like quantity
$\eta_R$ (the Legendre transform of $\pi_R$) introduced in \cite{H}. For
$R>0$ fixed let $\bigl(\cA_R^{(n)}\bigr)^{sa}$ and $(M_{N,R}^{sa})^n$ be as
in \S1. For each $h\in\bigl(\cA_R^{(n)}\bigr)^{sa}$  the {\it free pressure}
$\pi_R(h)$ of $h$ is defined by
$$
\pi_R(h):=\limsup_{N\to\infty}\biggl({1\over N^2}P_{N,R}(h)
+{n\over2}\log N\biggr),
$$
where the (microstates) pressure function $P_{N,R}(h)$ is given as
$$
P_{N,R}(h):=\log\int_{(M_{N,R}^{sa})^n}
\exp\bigl(-N\Tr_N(h(A_1,\dots,A_n))\bigr)
\,d\Lambda_N^{\otimes n}(A_1,\dots,A_n).
$$
Note that $\pi_R$ is a convex function on $\bigl(\cA_R^{(n)}\bigr)^{sa}$ such
that $|\pi_R(h_1)-\pi_R(h_2)|\le\|h_1-h_2\|_R$ for all
$h_1,h_2\in\bigl(\cA_R^{(n)}\bigr)^{sa}$. For
$\tau\in TS\bigl(\cA_R^{(n)}\bigr)$ the quantity $\eta_R(\tau)$ is defined by
$$
\eta_R(\tau):=\inf\Bigl\{\tau(h)+\pi_R(h):
h\in\bigl(\cA_R^{(n)}\bigr)^{sa}\Bigr\}.
$$
We then have
$$
\pi_R(h)=\max\Bigl\{-\tau(h)+\eta_R(\tau):
\tau\in TS\bigl(\cA_R^{(n)}\bigr)\Bigr\}
$$
so that $\pi_R$ on $\bigl(\cA_R^{(n)}\bigr)^{sa}$ and $\eta_R$ on
$TS\bigl(\cA_R^{(n)}\bigr)$ are the Legendre transforms of each other with
respect to the Banach space duality between $\bigl(\cA_R^{(n)}\bigr)^{sa}$
and $\bigl(\cA_R^{(n)}\bigr)^{*,\,sa}$ ($\supset TS\bigl(\cA_R^{(n)}\bigr)$).
We say that $\tau\in TS\bigl(\cA_R^{(n)}\bigr)$ is an {\it equilibrium
tracial state} associated with $h\in\bigl(\cA_R^{(n)}\bigr)^{sa}$ if the
equality
\begin{equation}\label{F-3.1}
\pi_R(h)=-\tau(h)+\eta_R(h)
\end{equation}
holds. This equality is a kind of variational principle.

In this section we will prove the next TCI for non-commutative
tracial distributions with $\eta_R$ in place of $\chi$. Since
$\chi_R(\tau)\le\eta_R(\tau)$ \cite[Theorem 4.5]{H}, this TCI is sharper
than that given in Theorem \ref{T-2.2} though $\tau$ becomes restrictive
here. But it is worth noting (see \cite{BroRoc}, also \cite[V.1.1]{Is}) that
the set of $\tau\in TS\bigl(\cA_R^{(n)}\bigr)$ satisfying the assumption in
the theorem is norm-dense in
$\bigl\{\tau\in TS\bigl(\cA_R^{(n)}\bigr):\eta_R(\tau)>-\infty\bigr\}$.

\begin{thm}\label{T-3.1}
Let $\bQ=(Q_1,\dots,Q_n)$ be an $n$-tuple of real-valued continuous functions
on $\bR$, and assume that there exists a constant $\rho>0$ such that all
$Q_i(x)-{\rho\over2}x^2$, $1\le i\le n$, are convex on $\bR$. Assume $R\ge R_0$
with $R_0$ given in \S2. If $\tau\in TS\bigl(\cA_R^{(n)}\bigr)$ is an
equilibrium tracial state associated with some
$h\in\bigl(\cA_R^{(n)}\bigr)^{sa}$, then
\begin{equation}\label{F-3.2}
W_2(\tau,\tau_\bQ)\le\sqrt{{2\over\rho}\Biggl(-\eta_R(\tau)
+\tau\Biggl(\sum_{i=1}^nQ_i(X_i)\Biggr)+B_\bQ\Biggr)},
\end{equation}
where $B_\bQ$ is the constant in \eqref{F-2.3}.
\end{thm}

The essence of the proof is same as that of Theorem \ref{T-2.2} based on the random matrix approximation
procedure. For each $n\in\bN$ and $h\in\bigl(\cA_R^{(n)}\bigr)^{sa}$ define
$\lambda_{N,R}(h)\in\cP\bigl((M_{N,R}^{sa})^n\bigr)$ by
$$
{d\lambda_{N,R}(h)\over d\Lambda_N^{\otimes n}}(A_1,\dots,A_n)
={1\over Z_{N,R}(h)}\exp\bigl(-N\Tr_N(h(A_1,\dots,a_n))\bigr)
\chi_{(M_{N,R}^{sa})^n}(A_1,\dots,A_n)
$$
with the normalization constant $Z_{N,R}(h):=\exp\bigl(P_{N,R}(h)\bigr)$.
This $\lambda_{N,R}(h)$ is a unique probability measure on
$(M_{N,R}^{sa})^n$ satisfying the (microstates) Gibbs variational principle
\begin{equation}\label{F-3.3}
P_{N,R}(h)=-N^2\widehat\lambda_{N,R}(h)(h)
+S(\lambda_{N,R}(h))
\end{equation}
with the Boltzmann-Gibbs entropy
$$
S(\lambda_{N,R}(h)):=-\int_{(M_{N,R}^{sa})^n}
{d\lambda_{N,R}(h)\over d\Lambda_N^{\otimes n}}
\log{d\lambda_{N,R}(h)\over d\Lambda_N^{\otimes n}}
\,d\Lambda_N^{\otimes n}.
$$

\bigskip\noindent
{\it Proof of Theorem 3.1.}\enspace
We may prove that
\begin{equation}\label{F-3.4}
{1\over2}W_2(\tau_0,\tau_\bQ)^2\le{1\over\rho}
\Biggl(-\pi_R(h_0)+\tau_0\Biggl(\sum_{i=1}^nQ_i(X_i)-h_0\Biggr)
+B_\bQ\Biggr)
\end{equation}
when $\tau_0\in TS\bigl(\cA_R^{(n)}\bigr)$ and
$h_0\in\bigl(\cA_R^{(n)}\bigr)^{sa}$ satisfy the variational equality
\eqref{F-3.1}. Let us first assume that $\tau_0$ is a unique equilibrium
tracial state associated with $h_0$ (equivalently, $\pi_R$ is differentiable
at $h_0$). Choose a subsequence $N(1)<N(2)<\cdots$ such that
\begin{equation}
\pi_R\left(h_0\right)=\lim_{k\to\infty}\label{F-3.5}
\biggl({1\over N(k)^2}P_{N(k),R}(h_0)+{n\over2}\log N(k)\biggr)
\end{equation}
and $\widehat\lambda_{N(k),R}(h_0)$ weakly* converges to some
$\tau_1\in TS\bigl(\cA_R^{(n)}\bigr)$. For every
$h\in\bigl(\cA_R^{(n)}\bigr)^{sa}$ we get
\begin{eqnarray*}
&&\widehat\lambda_{N(k),R}(h_0)(h)
+{1\over N(k)^2}P_{N(k),R}(h) \\
&&\qquad\ge{1\over N(k)^2}S\bigl(\lambda_{N(k),R}(h_0)\bigr)
=\widehat\lambda_{N(k),R}(h_0)(h_0)
+{1\over N(k)^2}P_{N(k),R}(h_0)
\end{eqnarray*}
thanks to \eqref{F-3.3}. From this and \eqref{F-3.5} as well as the weak*
convergence of $\widehat\lambda_{N(k),R}(h_0)$ it is easy to see that
$$
\tau_1(h)+\pi_R(h)\ge\tau_1(h_0)+\pi_R(h_0)
$$
so that $\tau_1$ is an equilibrium tracial state associated with $h_0$.
Therefore,
\begin{equation}\label{F-3.6}
\widehat\lambda_{N(k),R}(h_0)\longrightarrow\tau_0
\ \,\mbox{weakly*\quad as $k\to\infty$}.
\end{equation}

For every $N\in\bN$ let $\lambda_N(\bQ)$, $Z_N(\bQ)$ and
$\widehat\lambda_{N,R}(\bQ)$ be defined as in \S2. By Lemmas \ref{L-1.3} and
\ref{L-2.4}, as in the proof of Theorem \ref{T-2.2} we have
\begin{eqnarray*}
&&{1\over2}W_2\bigl(\widehat\lambda_{N,R}(h_0),
\widehat\lambda_{N,R}(\bQ)\bigr) \\
&&\qquad\le{1\over\rho N^2}\int_{\left(M_{N,R}^{sa}\right)^n}
\log{d\lambda_{N,R}(h_0)\over d\lambda_N(\bQ)}
\,d\lambda_{N,R}\left(h_0\right) \\
&&\qquad={1\over\rho}\Biggl(-{1\over N^2}
S(\lambda_{N,R}(h_0))
+\widehat\lambda_{N,R}(h_0)
\Biggl(\sum_{i=1}^nQ_i(X_i)\Biggr)
+{1\over N^2}\log Z_N(\bQ)\Biggr) \\
&&\qquad={1\over\rho}\Biggl(-{1\over N^2}P_{N,R}(h_0)
+\widehat\lambda_{N,R}(h_0)
\Biggl(\sum_{i=1}^nQ_i(X_i)-h_0\Biggr)
+{1\over N^2}\log Z_N(\bQ)\Biggr)
\end{eqnarray*}
thanks to \eqref{F-3.1}. Now, restrict the above estimates to
the subsequence $N(1)<N(2)<\cdots$ and apply \eqref{F-3.5} as well as
\eqref{F-2.3}. By $2^\circ$ in \S\S1.3 together with Lemma \ref{L-2.1} and
\eqref{F-3.6}, we then obtain
\begin{eqnarray*}
{1\over2}W_2\bigl(\tau_0,\tau_\bQ\bigr)^2
&\le&\liminf_{k\to\infty}{1\over2}
W_2\bigl(\widehat\lambda_{N(k),R}(h_0),
\widehat\lambda_{N(k),R}(\bQ)\bigr)^2 \\
&\le&-\pi_R(h_0)
+\tau_0\Biggl(\sum_{i=1}^nQ_i(X_i)-h_0\Biggr)
+B_\bQ.
\end{eqnarray*}

Next, assume that $\tau_0$ is a not necessarily unique equilibrium tracial
state associated with $h_0$. According to \cite{LR} (also \cite[6.2.43]{BR}),
$\tau_0$ belongs to the weakly* closed convex hull of the set $\mathcal{T}_0$
of $\tau\in TS\bigl(\cA_R^{(n)}\bigr)$ for which there exist
$h_k\in\bigl(\cA_R^{(n)}\bigr)^{sa}$ and $\tau_k\in TS\bigl(\cA_R^{(n)}\bigr)$
such that $\tau_k$ is a unique equilibrium tracial state associated with $h_k$
for each $k\in\bN$,
$\|h_k-h_0\|_R\to0$ and $\tau_k\to\tau$ weakly*. To show \eqref{F-3.4} for
$\tau_0$ and $h_0$, it suffices thanks to $2^\circ$ and $3^\circ$ in \S\S1.3
to prove it for every $\tau\in\mathcal{T}_0$ and $h_0$. Let $h_k$ and
$\tau_k$ be as in the description of the set $\mathcal{T}_0$. Then, the
above-proven case implies that
$$
{1\over2}W_2(\tau_k,\tau_\bQ)^2
\le-\pi_R(h_k)+\tau_k\Biggl(\sum_{i=1}^nQ_i(X_i)-h_k\Biggr)+B_\bQ
$$
for all $k\in\bN$. Hence \eqref{F-3.4} for $\tau$ and $h_0$ is
obtained  by letting $k\to\infty$ in view of $2^\circ$ in \S\S1.3 and the
norm-continuity of $\pi_R$; thus the proof is
completed.\qed

\begin{cor}
Let $\bQ=(Q_1,\dots,Q_n)$ be an $n$-tuple of real-valued continuous functions
on $\bR$ with the same assumption as in Theorem 3.1 for some $\rho>0$, and let
$R>0$ be as in Theorem 3.1. Then $\tau_\bQ$ is a unique equilibrium tracial
state associated with
$\sum_{i=1}^nQ_i(X_i)\in\bigl(\cA_R^{(n)}\bigr)^{sa}$.
\end{cor}

\proof
If $\tau_0\in TS\bigl(\cA_R^{(n)}\bigr)$ is an equilibrium tracial state
associated with $\sum_{i=1}^nQ_i(X_i)$, then the right-hand side of
\eqref{F-3.2} (or \eqref{F-3.4}) is zero so that $\tau_0=\tau_\bQ$.\qed

\bigskip
In particular, let $Q_i(x)=x^2/2$, $1\le i\le n$, and $R\ge2$. 
If $\tau\in TS\bigl(\cA_R^{(n)}\bigr)$ is an equilibrium tracial state
associated with some $h\in\bigl(\cA_R^{(n)}\bigr)^{sa}$, then 
$$
W_2\Bigl(\tau,\gamma_{0,2}^{\bigstar n}\Bigr)
\le \sqrt{2\Biggl(-\eta_R(\tau)
+\tau\Biggl({1\over2}\sum_{i=1}^n X_i^2\Biggr)
+{n\over2}\log2\pi\Biggr)},
$$
where $\gamma_{0,2}^{\bigstar n}$ denotes the distribution of a standard
semicircular system. Hence, $\gamma_{0,2}^{\bigstar n}$ is a unique
equilibrium tracial state associated with ${1\over2}\sum_{i=1}^nX_i^2$. This
also says that $\eta_R$ admits a unique maximizer $\gamma_{0,2}^{\bigstar n}$
when restricted on
$\bigl\{\tau\in TS\bigl(\cA_R^{(n)}\bigr):\tau\bigl(\sum_{i=1}^n
X_i^2\bigr)\le n\bigr\}$, which is a refinement of the same result for
$\chi$ in \cite{V4}.

The question whether the TCI \eqref{F-3.2} holds or not for any
$\tau\in TS\bigl(\cA_R^{(n)}\bigr)$ without the equilibrium assumption is 
still left open (and seems very important to obtain an in-depth
understanding of $\eta_R$).

\section{The unitary case}
\setcounter{equation}{0}

\subsection{TCI for $\chi_u$}
For each $n\in\bN$, the universal free product $C^*$-algebra
$C(\bT)^{\bigstar n}$, where $\bT$ is the unit circle, is nothing but the
universal group $C^*$-algebra $C^*(\bF_n)$ of the free group $\bF_n$ of $n$
generators. Let $g_1,\dots,g_n$ denote the canonical $n$ unitary generators
of $C^*(\bF_n)$. For each $\tau\in TS(C^*(\bF_n))$ take the $W^*$-probability
space $\bigl(\pi_\tau(C^*(\bF_n))'',\tilde\tau)\bigr)$ via the GNS
representation $\pi_\tau$ and define the {\it free entropy} (unitary version)
$\chi_u(\tau)$ by
$$
\chi_u(\tau):=\chi_u(\pi_\tau(g_1),\dots,\pi_\tau(g_n))
$$
(see \cite[\S6.5]{HP} for the precise definition of the microstates free
entropy for $n$-tuples of unitaries).

On the other hand, the $p$-Wasserstein distance
$W_p(\tau_1,\tau_2)$ between $\tau_1,\tau_2\in TS(C^*(\bF_n))$ is defined by
\eqref{F-1.2} with $(g_1,\dots,g_n)$ in place of $(a_1,\dots,a_n)$. (Note
that $a_1,\dots,a_n$ were not necessarily self-adjoint in \S\S1.3.)

For each real-valued continuous function $Q$ on $\bT$, the functional
$$
-\Sigma(\mu)+\mu(Q)
:=-\iint_{\!\!\!\bT^2}\log|\zeta-\eta|\,d\mu(\zeta)\,d\mu(\eta)
+\int_\bT Q(\zeta)\,d\mu(\zeta),\quad \mu\in\cP(\bT), 
$$
has a unique minimizer $\mu_Q$ called the {\it equilibrium measure} associated
with $Q$ (see \cite{ST}). When $\bQ=(Q_1,\dots,Q_n)$ is an $n$-tuple of
real-valued continuous functions on $\bT$, we define
$\tau_\bQ\in TS(C^*(\bF_n))$ as the free product of $\mu_{Q_i}$'s, i.e.,
$\tau_\bQ:=\bigstar_{i=1}^n\mu_{Q_i}$ on $C^*(\bF_n)=C(\bT)^{\bigstar n}$.

The next theorem is the counterpart of Theorem \ref{T-2.2} in the unitary
setting.

\begin{thm}\label{T-4.1}
Assume that there exists a constant $\rho>-{1\over2}$ such that all
$Q_i\bigl(e^{\im t}\bigr)-{\rho\over2}t^2$, $1\le i\le n$, are convex on
$\bR$. Then we have
$$
W_2(\tau,\tau_\bQ)\le\sqrt{{4\over1+2\rho}
\Biggl(-\chi_u(\tau)+\tau\Biggl(\sum_{i=1}^nQ_i(g_i)\Biggr)+B_\bQ\Biggr)}
$$
for every $\tau\in TS(C^*(\bF_n))$, where
$$
B_\bQ:=\chi_u(\tau_\bQ)-\tau_\bQ\Biggl(\sum_{i=1}^nQ_i(g_i)\Biggr)
$$
{\rm(}See also $3^\circ$ below for the constant $B_\bQ${\rm)}.
Furthermore, $\tau_\bQ$ is a unique minimizer of
$-\chi_u(\tau)+\tau\bigl(\sum_{i=1}^nQ_i(g_i)\bigr)$ for
$\tau\in TS(C^*(\bF_n))$. 
\end{thm}

In the special case where $Q_i$'s are all zero and so $\rho=0$, the above
inequality becomes
$$
W_2\Bigl(\tau,\gamma_0^{\bigstar n}\Bigr)\le2\sqrt{-\chi_u(\tau)},
\qquad\tau\in TS(C^*(\bF_n)),
$$
where the free product state $\gamma_0^{\bigstar n}$ is the distribution 
of a standard Haar unitary system of $n$ variables.

A key idea in proving the theorem is to apply the classical TCI in the
Riemannian setting in a certain random matrix approximation. Here, by a
geometric reason on Ricci curvature tensors, random matrices at our disposal
are special unitary ones instead of unitary. Some important facts
needed in the proof are in order.

\medskip
$1^\circ$ {\bf The $\SU$-microstates free entropy.}\enspace
Even when $\U(N)$ is replaced by $\SU(N)$ in the definition of
$\chi_u(u_1,\dots,u_n)$ \cite[\S6.5]{HP}, the microstates free entropy
introduced is the same. To prove this, define
$\xi:\bT^N\to\{(\zeta_1,\dots,\zeta_N)\in\bT^N:\zeta_1\cdots\zeta_N=1\}$ by
$$
\xi(\zeta_1,\dots,\zeta_N):=\bigl(\zeta_1(\zeta_1\cdots\zeta_N)^{-1/N},
\dots,\zeta_N(\zeta_1\cdots\zeta_N)^{-1/N}\bigr),
$$
where $\zeta^{1/N}$ for $\zeta\in\bT$ is the principal $N$th root and
$\zeta^{-1/N}:=(\zeta^{1/N})^{-1}$, and define $\Xi:\U(N)\to\SU(N)$ by
$\Xi(U):=V{\rm diag}\,\xi(\zeta_1,\dots,\zeta_N)V^*$ under a diagonalization
$U=V{\rm diag}(\zeta_1,\dots,\zeta_N)V^*$ with $V\in\U(N)$ and
$(\zeta_1,\dots,\zeta_N)\in\bT^N$. Then, $\Xi$ is a well-defined Borel
measurable map and we have $\gamma_N^\U\circ\Xi^{-1}=\gamma_N^\SU$ (see
\S\S1.1 for notations). Now, the above-mentioned fact can be directly shown
by using the forms of $\gamma_N^\U$ and $\gamma_N^\SU$ under diagonalizations
(see e.g.~\cite[\S\S1.5]{HPU2}).

\medskip
$2^\circ$ {\bf A key inequality.}\enspace
For each $\lambda\in\cP(\SU(N)^n)$ define the distribution
$\widehat\lambda\in TS(C^*(\bF_n))$ by
$$
\widehat\lambda(h)
:=\int_{\SU(N)^n}{1\over N}\Tr_N(h(U_1,\dots,U_n))
\,d\lambda(U_1,\dots,U_n), \quad h\in C^*(\bF_n) ,
$$
where $h\in C^*(\bF_n))\mapsto h(U_1,\dots,U_n)\in M_N(\bC)$ is the
$\ast$-homomorphism (``non-commutative functional calculus") sending each
$g_i$ to $U_i$ for each $(U_1,\dots,U_n)\in\SU(N)$. For every
$\lambda_1,\lambda_2\in\cP(\SU(N)^n)$ we have
$$
W_2\bigl(\widehat\lambda_1,\widehat\lambda_2\bigr)
\le{1\over\sqrt N}W_{2,\,{HS}}(\lambda_1,\lambda_2)
\le{1\over\sqrt N}W_{2,\,{\rm geod}}(\lambda_1,\lambda_2),
$$
where $W_{2,\,{HS}}(\lambda_1,\lambda_2)$ is the Wasserstein
distance with respect to the distance on $\SU(N)^n$ induced by the
Hilbert-Schmidt norm, while $W_{2,\,{\rm geod}}(\lambda_1,\lambda_2)$ with
respect to the geodesic distance. The proof of the first inequality is
similar to that of Lemma
\ref{L-1.3} while the second is obvious.

\medskip
$3^\circ$ {\bf Asymptotic freeness for $\SU$-random matrices.}\enspace
Let $\bQ=(Q_1,\dots,Q_n)$ be real-valued continuous functions on $\bT$. For
each
$n\in\bN$ define $\lambda_N(\bQ)\in\cP(\SU(N)^n)$ by the product measure
$\lambda_N(\bQ):=\bigotimes_{i=1}^n\lambda_N(Q_i)$ of
$$
d\lambda_N(Q_i)(U):={1\over Z_N(Q_i)}
\exp\bigl(-N\Tr_N(Q_i(U))\bigr)\,d\gamma_N^\SU(U)
$$
with a normalization constant $Z_N(Q_i)$. The asymptotic freeness for
unitary random matrices due to \cite{Vo} remains valid for special unitary
random matrices. In fact, a stronger result on the almost sure asymptotic
freeness for independent special unitary random matrices can be shown by
modifying the proof in \cite[4.3.5]{HP}. Also, as a consequence of the
large deviation theorem \cite[Theorem 2.1]{HPU3}, it follows that
the mean eigenvalue distribution of $\lambda_N(Q_i)$ converges
to $\mu_{Q_i}$ for $1\le i\le n$. We thus see that
$$
\widehat\lambda_N(\bQ)\longrightarrow\tau_\bQ=\bigstar_{i=1}^n\mu_{Q_i}
\ \ \mbox{weakly*}.
$$
Moreover, since
$$
\lim_{N\to\infty}{1\over N^2}\log Z_N(Q_i)
=\Sigma(\mu_{Q_i})-\mu_{Q_i}(Q_i)
$$
(see \cite[Theorem 2.1]{HPU3}), we notice that the constant
$B_\bQ$ in Theorem \ref{T-4.1} can be expressed as 
$$
B_\bQ=\sum_{i=1}^n\bigl(\Sigma(\mu_{Q_i})-\mu_{Q_i}(Q_i)\bigr)
=\lim_{N\to\infty}{1\over N^2}\sum_{i=1}^n\log Z_N(Q_i).
$$

\medskip
$4^\circ$ {\bf TCI on $\SU(N)^n$.}\enspace
Let $\bQ=(Q_1,\dots,Q_n)$ be as in Theorem \ref{T-4.1} and assume further
that all $Q_i$'s are $C^2$-functions. Then, thanks to \cite[Lemmas 1.2 and
1.3]{HPU2}, the function
$$
\Psi_N(U_1,\dots,U_n):=N\Tr_N\Biggl(\sum_{i=1}^nQ_i(U_i)\Biggr)
$$
is a $C^2$-function on $\SU(N)^n$ with the Hessian
$\Hess(\Psi_N)\ge N\rho I_{(N^2-1)n}$. Also, note that the Ricci
curvature tensor of $\SU(N)^n$ is $\Ric(\SU(N)^n)={N\over2}I_{(N^2-1)n}$.
Hence, by the TCI in the Riemannian manifold setting due to \cite{OV} combined
with \cite{BE}, we obtain
$$
W_{2,\,{\rm geod}}(\lambda,\lambda_N(\bQ))
\le\sqrt{{4\over N(1+2\rho)}S(\lambda,\lambda_N(\bQ))}
$$
for every $\lambda\in\cP(\SU(N)^n)$.

\medskip
Now, the proof of Theorem \ref{T-4.1} based on the above facts
$1^\circ$--$4^\circ$ is analogous to that of Theorem \ref{T-2.2}; so the
details are left to the reader. But, it is worthwhile to note one more point.
As in the proof of \cite[Theorem 2.7]{HPU2}, the regularization technique by
the  use of Poisson integrals enables us to assume that all $Q_i$'s are
smooth  functions on $\bT$ so that one can go through with $4^\circ$. 

\subsection{TCI for $\eta_u$}

For each $h\in C^*(\bF_n)^{sa}$ we introduce the {\it free pressure} (unitary
version) $\pi_u(h)$ by
\begin{eqnarray*}
&&\pi_u(h) \\
&&\quad:=\limsup_{N\to\infty}\Biggl({1\over N^2}\log
\int_{\U(N)^n}\exp\bigl(-N\Tr_N(h(U_1,\dots,U_n))\bigr)
\,d(\gamma_N^\U)^{\otimes n}(U_1,\dots,U_n)\Biggr) \\
&&\quad\ =\limsup_{N\to\infty}\Biggl({1\over N^2}\log
\int_{\SU(N)^n}\exp\bigl(-N\Tr_N(h(U_1,\dots,U_n))\bigr)
\,d(\gamma_N^\SU)^{\otimes n}(U_1,\dots,U_n)\Biggr).
\end{eqnarray*}
The equality of the two $\limsup$'s can be shown from the fact stated in the
above $1^\circ$. As in the self-adjoint setting \cite[Proposition 2.3]{H},
$\pi_u$ is convex on $C^*(\bF_n)^{sa}$ and
$|\pi_u(h_1)-\pi_u(h_2)|\le\|h_1-h_2\|$ for all $h_1,h_2\in C^*(\bF_n)^{sa}$.
It is seen as in \cite[Theorem 3.4]{H} (or \eqref{F-3.1}) that $\pi_u$ is the
converse Legendre transform of $\eta_u$ as
$$
\pi_u(h)=\max\bigl\{-\tau(h)+\eta_u(\tau):
\tau\in TS(C^*(\bF_n))\bigr\},
\qquad h\in C^*(\bF_n)^{sa},
$$
and we say that $\tau$ is an {\it equilibrium tracial state} associated with
$h$ if $\pi_u(h)=-\tau(h)+\eta_u(h)$ holds.

In particular when $N=1$ and $\mu\in\cP(\bT)$ ($=TS(C(\bT))$) we have
$\eta_u(\mu)=\chi_u(\mu)$ ($=\Sigma(\mu)$) (see \cite[\S6]{HMP}). The proof
of the next theorem is similar to that of \cite[Theorem 4.5]{H}.

\begin{thm}\label{T-4.2}
We have $\chi_u(\tau)\le\eta_u(\tau)$ for every $\tau\in TS(C^*(\bF_n))$.
Moreover, if $\tau$ is a free product tracial state {\rm (}i.e.,
$g_1,\dots,g_n$ are $\ast$-free with respect to $\tau${\rm )}, then
$\chi_u(\tau)=\eta_u(\tau)$.
\end{thm}

Furthermore, by considering the minimal $C^*$-tensor product
$C^*(\bF_n)\otimes_{\rm min}C^*(\bF_n)$ as in \cite[\S6]{H}, the definition
of $\eta_u$ can be modified so that the modified $\tilde\eta_u(\tau)$ is
equal to $\chi_u(\tau)$ for all $\tau\in TS(C^*(\bF_n))$. This result is of
considerable importance but it is not directly related to the free TCI in
Theorem \ref{T-4.3}; so we omit the details.

Finally, we state the counterpart of Theorem \ref{T-3.1} in the unitary
setting; the TCI is sharper than that in Theorem \ref{T-4.1} though $\tau$ is
rather restricted. The structure of the proof is quite parallel with that of
Theorem \ref{T-3.1} and the details are again left to the reader.

\begin{thm}\label{T-4.3}
Let $\bQ=(Q_1,\dots,Q_n)$ be real-valued continuous functions on $\bT$
satisfying the same assumption as in Theorem \ref{T-4.1} with a
constant $\rho>-{1\over2}$. If $\tau\in TS(C^*(\bF_n))$ is an equilibrium
tracial state associated with some $h\in C^*(\bF_n)^{sa}$, then
$$
W_2(\tau,\tau_\bQ)\le\sqrt{{4\over1+2\rho}
\Biggl(-\eta_u(\tau)+\tau\Biggl(\sum_{i=1}^n Q_i(g_i)\Biggr)+B_\bQ\Biggr)},
$$
where $B_\bQ$ is the same constant as in Theorem \ref{T-4.1}.
Furthermore, $\tau_\bQ$ is a unique equilibrium tracial state associated with
$\sum_{i=1}^n Q_i(g_i)\in C^*(\bF_n)^{sa}$.
\end{thm}

\end{document}